\documentclass[11pt]{article}
\usepackage{amssymb}
\usepackage{mathrsfs}
\usepackage{indentfirst}
\usepackage{latexsym,amsmath,amssymb,amsfonts,epsfig,graphicx,cite,psfrag}
\usepackage{eepic,color,colordvi,amscd}

\newtheorem{theo}{Theorem}[section]

\newtheorem{prop}[theo]{Proposition}
\newtheorem{conj}[theo]{Conjecture}

\begin{document}

\title{Monochromatic loose path partitions in $k$-uniform hypergraphs}
\author{Changhong Lu, Bing Wang, Ping Zhang\\
Department of Mathematics\\
East China Normal University\\
Shanghai 200241, China\\
}
\date{}

\maketitle

\begin{abstract}
A conjecture of Gy\'{a}rf\'{a}s and S\'{a}rk\"{o}zy says that in every
$2$-coloring of the edges of the complete $k$-uniform hypergraph
$K_n^k$, there are two disjoint monochromatic loose paths of
distinct colors such that they cover all but at most $k-2$ vertices.
A weaker form of this conjecture with $2k-5$ uncovered vertices instead of
$k-2$ is proved, thus the conjecture holds for $k=3$. The main
result of this paper states that the conjecture is true for all $k\ge 3$.

\bigskip

\end{abstract}

\footnote{E-mail addresses: chlu@math.ecnu.edu.cn (C. Lu), wuyuwuyou@126.com (B. Wang), mathzhangping@126.com (P. Zhang).}

\section{Introduction}

A hypergraph $H=(V,E)$ consists of a set $V$ of vertices and a set
$E$ of edges, where each edge is a subset of $V$. If all the edges
of $H$ have same size $k$, then the hypergraph $H$
is said to be $k$-uniform. Let $K_n^k$ denote the complete $k$-uniform
hypergraph on $n$ vertices (the family of all $k$-element subsets of
a $n$-element set). A $k$-uniform loose (or linear) path of length $\ell$,
denoted $\mathcal{P}_{\ell}^k$, is a $k$-uniform hypergraph with edges
$e_1,e_2,\cdots,e_\ell$ such that $\forall i\in [\ell-1], |e_i\cap e_{i+1}|=1$
and $|e_i\cap e_j|=0$ for all other pairs $\{i,j\}$, $i\ne j$. For a loose path
$\mathcal{P}_{\ell}^k$ and a vertex $v\in V(\mathcal{P}_{\ell}^k)$, if $v$ lies
in two edges of $\mathcal{P}_{\ell}^k$, then we call $v$ a $2$-degree vertex of
$\mathcal{P}_{\ell}^k$. A $k$-uniform tight path of length $\ell$, is a
sequence of $k+\ell$ vertices with every
consecutive set of $k$ vertices forms an edge. For $k=2$ we
obtain the usual definition of a path $P_\ell$ with $\ell$ edges.

In this paper $r$-coloring always means edge-coloring with $r$ colors
(traditionally red and blue when $r=2$). The following simple proposition,
introduced by Gerenc\'{s}er and Gy\'{a}rf\'{a}s in \cite{LA1},
is our starting point here.

\begin{prop}\label{prop1}
In any 2-coloring of the edges of a finite complete graph the vertices can
be partitioned into a red and a blue path. Here the empty graph and the
one-vertex graph is accepted as a path of any color.
\end{prop}

Note that any result about covering the vertices of edge-colored graphs
 by a small number of monochromatic subgraphs will imply a Ramsey-type
 result as a corollary. For example, Proposition \ref{prop1}
implies the bound $R(P_n,P_m)\le n+m-3$ for $n,m\ge 2$. In fact, Proposition \ref{prop1}
subsequently gave birth to the area of partitioning edge-colored complete graphs
into monochromatic subgraphs. There have been many further results, questions and conjectures in
this area, many of which generalize Proposition \ref{prop1} in graphs or hypergraphs.
we refer to two surveys \cite{LA2,KL}. However, in contrast to
the graph case, there are only a few results on covering the vertices with
monochromatic pieces of hypergraphs, see for example, \cite{AG,G,LA3,VM}.

There are various definitions of paths and cycles (for example, Berge, loose and tight)
for hypergraphs. We focus on loose path here. Similar to the graph case,
a set of less than $k$ vertices in an edge-colored $k$-uniform hypergraph is
accepted as a loose path of any color. However, it seems difficult to extend
Proposition \ref{prop1} to loose or tight paths of hypergraphs. The following
conjecture first presented by Gy\'{a}rf\'{a}s and S\'{a}rk\"{o}zy can be found
in \cite{LA2} and \cite{LA3}.

\begin{conj}\label{conj1}
In every $2$-coloring of the edges of $K_n^k$ there are two disjoint
monochromatic loose paths of distinct colors covering all but at most
$k-2$ vertices. This estimate is sharp for sufficiently
large $n$.
\end{conj}

Gy\'{a}rf\'{a}s and S\'{a}rk\"{o}zy \cite{LA3} presented the following
construction to show that if Conjecture \ref{conj1} holds, then it is best
possible for $n$ large enough: Consider the complete $k$-uniform hypergraph
$K_n^k$ with vertex bipartitions $Q$ and $S$, where $|Q|=(k-1)m+1$, $|S|=2(k-1)$
and $m\ge4(k-1)$. Then color all $k$-element subsets of $Q$ red and all uncolored
$k$-element subsets of $Q\cup S$ blue.

\section{Partitions by monochromatic loose paths}

In this section, we will prove the following slightly stronger result than Conjecture \ref{conj1}.

\begin{theo}\label{theo1}
Suppose that the edges of the complete $k$-uniform hypergraph
$K_n^k$ are colored with two colors, where $n\equiv2\bmod (k-1)$. Then $V(K_n^k)$
can be partitioned into two monochromatic loose paths of distinct colors.
\end{theo}

It is obvious that Theorem \ref{theo1} implies Conjecture \ref{conj1}: For each
$n\ne 2\bmod (k-1)$, removing at most $k-2$ vertices from $K_n^k$ will obtain
a smaller $K_{n'}^k$ with $n'\equiv2\bmod (k-1)$, then by Theorem \ref{theo1},
$V(K_{n'}^k)$ can be partitioned into two monochromatic loose paths of distinct colors.
That is, there are two disjoint monochromatic loose paths of distinct colors
such that they cover all but at most $k-2$ vertices of $K_n^k$.

{\bf Proof of Theorem \ref{theo1}.} Suppose the assertion is false. Then take
vertex disjoint red and blue loose paths $P_R$ and $P_B$ such that
they cover as many vertices as possible, and subject to this, the
difference between $|V(P_R)|$ and $|V(P_B)|$ is maximal.
Let $W$ be the set of vertices uncovered by the paths $P_R$ and $P_B$.
Without loss of generality suppose that $|V(P_R)|\ge|V(P_B)|$. Then
we have the following claim:

{\bf Claim} $|V(P_B)|=r(k-1)+1$ for some integer $r\ge 1$, that is, $P_B$ is proper.

If $P_B$ is not proper, then $|V(P_B)|\le k-1$. Note that now the red path
$P_R$ is proper and $W$ is not empty.
Then $|V(P_B)|+|W|=(k-1)s+1$ for some integer $s\ge 1$.
Since $|V(P_R)|+|V(P_B)|$ is maximal then $|V(P_B)|=k-1$.
Let $e=\{v_1,\cdots,v_k\}$ be the last edge with an $2$-degree
vertex $v_1$ of $P_R$. Let $\{u_1,\cdots,u_{k-1}\}$ be the
vertex set of $P_B$. Then we have $|W|=1$. Otherwise, $|W|=(s-1)(k-1)+1$
for some integer $s\ge 2$. Let $w_1,\cdots,w_k$ be $k$ vertices of $W$.
Then both edges $\{v_k,u_1,\cdots,u_{k-1}\}$ and $\{v_k,w_1,\cdots,w_k\}$
are blue, hence the two edges form a new blue path, say $P'_B$. Then $P_R-e$ and $P'_B$
cover more vertices, a contradiction. Let $w$ be the unique vertex of $W$.
Now we consider two cases as follows.

{\bf Case 1} $|V(P_R)|=k$, that is, $P_R$ is induced by an edge.

Then $\{v_1,\cdots, v_k$\} is the unique edge of $P_R$. It is easy to check that
for each $i\in [k]$, $\{v_i,u_1,\cdots,u_{k-1}\}$ is blue. Then
$\{v_2,\cdots,v_{k},w\}$ is red. Otherwise the two edges form a new blue
path covering all vertices, a contradiction.
Now a blue edge $\{v_1,u_1,\cdots,u_{k-1}\}$ and a red edge
$\{v_2,\cdots,v_{k},w\}$ cover all vertices of $K_n^k$, this contradicts the hypothesis.

{\bf Case 2} $|V(P_R)|=t(k-1)+1$ for some integer $t\ge2$, that is, $P_R$
contains at least two edges.

Let $f=\{x_1,\cdots,x_k\}$ be the first edge with a $2$-degree vertex $x_k$
of $P_R$. Note that $f_1=\{u_1,\cdots, u_{k-1},w\}$ is red and
$f_2=\{v_k,u_1,\cdots,u_{k-1}\}$ is blue. Then $f_3=\{x_1,v_3,\cdots, v_{k},w\}$ is blue too.
Otherwise, the red path $P_R-e+f_3+f_1$
together with the blue path $\{v_2\}$ cover all vertices of $K_n^k$.
By symmetry $f_4=\{v_2,x_2,\cdots, x_{k-1},w\}$ is blue. Now
three edges $f_2,f_3$ and $f_4$ induce a blue path. The blue path together with the red path
$P_R-e-f$ can cover all vertices of $K_n^k$. A contradiction. This completes the
proof of the claim.

The claim means that both two paths are proper. Then $|W|=s(k-1)$ for some
$s\ge1$, since $n\equiv2\bmod (k-1)$. Let $w_1,\cdots,w_{k-1}$ be $k-1$
vertices of $W$. We first show that $P_R$ contains at least two edges.
Otherwise, $P_R$ induced by an edge $\{v_1,\cdots,v_k\}$.
Let $\{u_1,\cdots, u_k\}$ be the unique edge of $P_B$.
Then similar to above, both edges $\{v_k,u_2,\cdots,u_k\}$
and $\{v_k,w_1,\cdots,w_{k-1}\}$ are blue and hence form a
blue path. The blue path together with the red path $\{v_1,\cdots,v_{k-1}\}$
will cover more vertices, a contradiction.



Let $f=\{x_1,\cdots,x_{k-1},x_k\}$ and $e=\{v_1,v_2,\cdots,v_k\}$ be the
first and last edges of $P_R$ respectively, where $x_k$ and $v_1$ are two
$2$-degree vertices of $P_R$ ($x_k=v_1$ is allowed).
Let $g=\{u_1,u_2,\cdots,u_k\}$ be the last edge of $P_B$. If $P_B$ is of length
at least two, then $u_1$ is a $2$-degree vertex of $P_B$.

 For convenience, let $X=\{x_1,\cdots,x_{k-1}\}$, $V=\{v_2,\cdots,v_k\}$,
 $U=\{u_2,\cdots,u_k\}$ and $W'=\{w_1,\cdots,w_{k-1}\}$.
For each element $Y\in \{X, V, U,W'\}$, let $Y_{i}$ denote an $i$-element
subset of $Y$. Specially, let $Y_{0}=\emptyset$. By the assumption of $P_R$ and $P_B$ , 
the following results are easy to check:

(i) for $i\in [k]\setminus\{1\}$, $\{v_i\}\cup W'$ is blue, $\{u_i\}\cup W'$ is
red; 

(ii) for $i,j\in [k]\setminus\{1\}$, $\{v_i,u_j\}\cup W'_{k-2}$ is blue;
by symmetry, for $i\in [k-1], j\in [k]\setminus\{1\}$, $\{x_i,u_j\}\cup W'_{k-2}$ is also blue;


(iii) for $i\in [k-1]$, $\{w_i\}\cup V$ is red; (otherwise, the blue path
$P_B+\{u_k,v_k\}\cup W'_{k-2}+\{w_i\}\cup V$ together with the red path 
$P_R-e$ will cover more vertices,
where $w_i\notin W'_{k-2}$)

(iv) for $i\in [k-1]$, $j\in [k]\setminus\{1\}$, $\{x_i,v_j\}\cup W'_{k-2}$
is blue (otherwise, $P_R-e
+\{x_i,v_j\}\cup W'_{k-2}+\{w_l\}\cup V$ is a longer red path,
where $w_l\notin W'_{k-2}$);

(v) for $i\in [k]\setminus\{1\}$, $j\in [k-1]$, $\{v_i ,w_{j}\}\cup X_{k-2}$
is blue (otherwise, the red path $P'_R=P_R-f
+\{v_i ,w_{j}\}\cup X_{k-2}+\{u_k\}\cup W'$ and the blue path $P'_B=P_B-g$ is 
a new covering with $|V(P'_R)|+|V(P'_B)|=|V(P_R)|+|V(P_B)|$ and 
$|V(P'_R)|-|V(P'_B)|>|V(P_R)|-|V(P_B)|$, a contradiction);
by symmetry, for $i,j\in [k-1]$, $\{x_i ,w_{j}\}\cup V_{k-2}$ is also blue;

(vi) for $i,j\in [k]\setminus\{1\}$, $l\in [k-1]$, $\{v_i,u_j,w_{l}\}\cup X_{k-3}$
is blue (otherwise, the red path $P_R-e
+\{v_i ,u_j,w_{l}\}\cup X_{k-3}+\{x_s\}\cup U$ and the blue path $P_B-g$ is a new covering, 
where $x_s\notin X_{k-3}$. Similar to (v), we get a contradiction).

Then $P_R'=P_R-f-e$ and
$P_B'=P_B+\{u_2,v_2,w_1\}\cup X_{k-3}+
\{w_1,x_i\}\cup V_{k-2}+\{v_j,x_l\}\cup W'_{k-2}$ can cover
more vertices than $P_R$ and $P_B$, where
$x_i,x_l\notin X_{k-3}$, $w_1\notin W'_{k-2}$ and
$v_j\in V_{k-2}$, a contradiction.
This completes the proof.

\end{document}